\documentclass[11pt, a4paper]{amsart}

\usepackage[english]{babel}
\usepackage{amsmath,enumerate, amsfonts, amssymb,amsthm, xypic}
\usepackage[dvips]{graphics} 

\newtheorem{thm}{Theorem}[section]
\newtheorem{lemma}[thm]{Lemma}
\newtheorem{prop}[thm]{Proposition}
\newtheorem{cor}[thm]{Corollary}

\theoremstyle{definition}
\newtheorem{defi}[thm]{Definition}

\newtheorem{rmk}[thm]{Remark}

\DeclareMathOperator{\jac}{jac}

\title{The corank is invariant under Blow-Nash equivalence}

\author{Goulwen Fichou}

\address {D\'epartement de Math\'ematiques,
Universit\'e d'Angers, 2 bd Lavoisier, 49045 Angers Cedex, France}

\email{fichou@tonton.univ-angers.fr}
\begin{document}

\begin{abstract} We address the following question, raised by T. Fukui.
  Is the corank an invariant of the blow-analytic equivalence between
  real analytic function germs? We give a partial positive answer in
  the particular case of the blow-Nash equivalence. The proof is based
  on the computation of some virtual Poincar\'e polynomials and zeta
  functions associated to a Nash function germ.
 \end{abstract}

\maketitle
\section*{Introduction}
The classification of real analytic function germs is a difficult
topic, notably in the choice of a good equivalence relation, between
germs, to study. An interesting relation, called blow-analytic equivalence, has been introduced by T. C. Kuo
\cite{kuo} and studied by several authors (see \cite{fukui} for a
recent survey). Notably, it has been proved that such an equivalence
relation does not admit moduli for a family with isolated
singularities. Moreover, the proof of this result produces effective
methods to prove blow-analytic triviality. On the other hand, some
invariants have been introduced in order
to distinguish blow-analytic types. However, because of the
complexity of these invariants, it remains difficult to obtain effective
classification results, at least in dimension greater than 3.

In this paper, we address the following related question, raised by
T. Fukui.

 Let $f:(\mathbb R^d,0) \longrightarrow (\mathbb R,0)$ be an
analytic germ. Assume that $0$ is singular for $f$, which means that
the jacobian matrix of $f$ at $0$ does vanish. Let $r$ denote the rank
of the hessian of $F$ at $0$. Then $f$ is
analytically equivalent to a function of the form
$$ \sum_{i=1}^s x_i^2-\sum_{j=s+1}^{s+t} x_j^2 + F (x),$$
where $s+t=r$ (note that $s$ or $t$ may vanish) and the order of $F$ is at least equal to $3$. The corank of $f$
is defined to be the corank of its hessian matrix at $0$, that is $d-r$.

\vskip 5mm

\textbf{Question:} Is the corank of an analytic function germ an
invariant of the blow-analytic equivalence?

\vskip 5mm

The answer to such a question would be a step toward a better
understanding of the blow-analytic equivalence relation, and therefore
to a better understanding of the singularities of real analytic
function germs.

We will not give a complete answer to this question, but we concentrate
in a particular case of the blow-analytic equivalence, where we
are able to conclude. Actually, a similar
situation holds in the Nash setting (that is analytic and moreover semi-algebraic). One can define the
blow-Nash equivalence between Nash function germs, and this relation
still has good triviality properties and effective invariants called zeta
functions \cite{fichou}. Recall that these zeta functions are defined
using an additive and multiplicative invariant (invariant means under Nash isomorphisms) of real algebraic sets,
the virtual Poincar\'e polynomial (cf. part \ref{zeta}).

The main result of this paper states that blow-Nash equivalent
Nash function germs have the same corank. Moreover, they have the same
index (the index corresponds to the integer $t$ above).

The proof is based on the invariance of the zeta functions
with respect to the blow-Nash equivalence \cite{fichou}, and on the
computation on a significant
part of these zeta functions for germs of Nash functions of the type
$\sum_{i=1}^s x_i^2-\sum_{j=s+1}^{s+t} x_j^2 $ (cf part
\ref{partcalculzeta}). To reach this aim, we need to compute some virtual
Poincar\'e polynomials associated to these germs. Such a computation
may be difficult in general, but here we manage to conclude thanks to the
degeneracy of a Leray-Serre spectral sequence (cf part \ref{partSS}).

\vskip 5mm

\textbf{Acknowledgements.} The author does wish to thank Toshizumi Fukui
for all his help concerning this work, work which has been initiated during
the invited stay of the author at Saitama University. He thanks also Adam
Parusi\'nski for valuable discussions.
\section{Corank and Blow-Nash equivalence}

\subsection{Blow-Nash equivalence}\label{zeta}

In this section, we recall briefly the notion of blow-Nash equivalence
as well as that of zeta functions. For more details, the reader is
refer to \cite{fichou}.

\begin{defi}\label{defiblownash}\begin{flushleft}\end{flushleft}
\begin{enumerate}
\item An algebraic modification of a
Nash function germ $f:(\mathbb R ^d,0)
\longrightarrow (\mathbb R,0)$ is a proper birational algebraic
morphism $\sigma_f:\big(M_f,\sigma_f^{-1}(0)\big)  \longrightarrow
(\mathbb R ^d,0) $, between Nash neighbourhoods of $0$ in $\mathbb
R ^d$ and the exceptional divisor $\sigma_f^{-1}(0)$ in $M_f$, which is
an isomorphism over the complement of the zero locus of $f$ and for which
 $f\circ \sigma$ is in normal crossing.
\item Let $f,g:(\mathbb R ^d,0) \longrightarrow (\mathbb R,0)$ be germs
of Nash functions. They are said to be blow-Nash equivalent if there
exist two algebraic modifications $$\sigma_f~:~\big(M_f,\sigma_f^{-1}(0)\big)  \longrightarrow
(\mathbb R ^d,0) \textrm{ and } \sigma_g:\big(M_g,\sigma_g^{-1}(0)\big)  \longrightarrow
(\mathbb R ^d,0),$$ such that $f\circ \sigma_f$ and $\jac \sigma_f$
(respectively $g\circ \sigma_g$ and $\jac \sigma_g$) have only normal
crossings simultaneously and a Nash isomorphism
(i.e. a semi-algebraic map which is an analytic isomorphism) $\Phi$
between analytic neighbourhoods $\big(M_f,\sigma_f^{-1}(0)\big)$ and $\big(M_g,\sigma_g^{-1}(0)\big)$ which
preserves the multiplicities of the jacobian determinants of $\sigma_f$ and
$\sigma_g$ along the components of the exceptional divisor, and which
induces a homeomorphism $\phi :(\mathbb R ^d,0) \longrightarrow
(\mathbb R^d,0)$ such that $f=g \circ \phi$, as illustrated by the commutative diagram:
$$\xymatrix{\big(M_f,\sigma_f^{-1}(0)\big) \ar[rr]^{\Phi} \ar[d]_{\sigma_f}& &\big(M_g,\sigma_g^{-1}(0)\big) \ar[d]^{\sigma_g}\\
            (\mathbb R ^d,0)       \ar[rr]^{\phi} \ar[dr]_f        & &           (\mathbb R ^d,0) \ar[dl]^g \\
         &   (\mathbb R,0) & } $$
\end{enumerate}
\end{defi}

The main results, concerning this equivalence relation between Nash
function germs, are that, on one hand, it does not admit moduli for
a Nash family with isolated singularity, and one the other hand, we
know invariants, called zeta functions, that will be crucial in the proof of the main result
of this paper. We recall now the definition of these
zeta functions. To begin with, let
us introduce the virtual Poincar\'e polynomial.

By an additive map on the category of real algebraic sets, we mean a
map $\beta$ such that $\beta(X)=\beta(Y) + \beta(X \setminus Y)$ where
$Y $ is an algebraic subset closed in $X$. Moreover $\beta$ is called
multiplicative if $\beta(X_1 \times X_2)=\beta(X_1) \cdot \beta(X_2)$
for real algebraic sets $X_1,X_2$.

\begin{prop}(\cite{fichou})\label{beta-nash} There exist additive maps
  on the category of real algebraic sets with values
  in $\mathbb Z$, denoted $\beta_i$ and called
  virtual Betti numbers, that coincide with the classical
  Betti numbers $\dim H_i(\cdot, \frac {\mathbb Z}{2 \mathbb Z})$ on
  the connected component of compact nonsingular real algebraic varieties.

Moreover $\beta(\cdot)=\sum_{i \geq 0} \beta_i(\cdot)u^i$ is
multiplicative, with values in $\mathbb Z [u]$.

Finally, if $X_1$ and $X_2$ are
   Nash isomorphic real algebraic sets, then $\beta(X_1)=\beta(X_2)$.
\end{prop}

Then we can define the zeta functions of a Nash function germ
$f:(\mathbb R ^d,0) \longrightarrow (\mathbb R,0)$ as follows. 
Denote by $\mathcal L$ the space of arcs at the
origin $0 \in \mathbb R ^d$, that is:
$$\mathcal L=\mathcal L(\mathbb R ^d,0)= \{\gamma : (\mathbb R,0) \longrightarrow (\mathbb R ^d,0)
:\gamma \textrm{ formal}\},$$
and by $\mathcal L_n$ the space of arcs truncated at the order $n+1$:
$$\mathcal L_n=\mathcal L_n(\mathbb R ^d,0)= \{\gamma \in \mathcal L:
\gamma (t)=a_1t+a_2t^2+ \cdots a_nt^n,~a_i \in \mathbb R ^d\},$$
for $n\geq 0$ an integer.

 We define the naive zeta function $Z_f(T)$ of $f$ as the
following element of $\mathbb Z[u,u^{-1}][[T]]$: 
$$Z_f(T)= \sum _{n \geq 1}{\beta (A_n)u^{-nd}T^n},$$
where 
$$A_n =\{\gamma \in  \mathcal L_n: ord(f\circ \gamma) =n \}=\{\gamma \in  \mathcal L_n: f\circ \gamma (t)=bt^n+\cdots,
b\neq 0 \}.$$
Similarly, we define zeta functions with sign
by
$$Z_f^{+1}(T)= \sum _{n \geq 1}{\beta (A_n^{+1})u^{-nd}T^n}
\textrm{~~and~~ }Z_f^{-1}(T)= \sum _{n \geq 1}{\beta (A_n^{-1})u^{-nd}T^n},$$
where 
$$A_n^{+1} =\{\gamma \in  \mathcal L_n: f\circ \gamma
(t)=+t^n+\cdots \}\textrm{~~and~~ } A_n^{-1} =\{\gamma \in  \mathcal L_n: f\circ \gamma
(t)=-t^n+\cdots \} .$$

The main result concerning these zeta functions is the following:

\begin{thm}(\cite{fichou})\label{inv} Blow-Nash equivalent Nash function germs
  have the same naive zeta function and the same zeta functions with sign. 
\end{thm}
\subsection{Corank of a Nash function germ}

The corank and the index of a Nash function germ is defined in a similar way than in
the analytic case. Moreover, the same splitting lemma holds in the
Nash case.

\begin{lemma}\label{morse} Let $f:(\mathbb R^d,0) \longrightarrow (\mathbb R,0)$
  be a Nash function germ. Assume that the jacobian matrix of $f$ does
  vanish at $0$. Then there exist a
  Nash isomorphism $\phi:(\mathbb R^d,0) \longrightarrow (\mathbb
  R^d,0)$ and integers $s,t$ (possibly equal to zero), where $s+t$ equals the
  rank of the hessian matrix of $f$ at $0$, such that
$$f \circ \phi(x_1,\ldots,x_d)= \sum_{i=1}^s x_i^2-\sum_{j=s+1}^{s+t}
x_j^2 + F (x_{s+t+1},\ldots, x_d),$$
where $F$ is a Nash function germ with order at least equal to $3$.
\end{lemma}

Classical proofs of lemma \ref{morse}, in the smooth case, use a
parametrized version of the Morse Lemma \cite{arnold}, but this is not allowed in
the Nash setting since it requires integration along vector
fields. However, an elementary proof, using only the implicit function
theorem and the Hadamard division lemma, has been given in \cite{LM}. This method adapts to our case since the implicit function
theorem does hold in the Nash \cite{BCR}, whereas the Hadamard division lemma is no longer necessary because
we are working with analytic functions.

%
%

Let us state now the central result of this paper.

\begin{thm}\label{mainthm} Let $f,g :(\mathbb R^d,0) \longrightarrow (\mathbb R,0)$
  be blow-Nash equivalent Nash function germs. Then $f$ and $g$ have
  the same corank and the same index.
\end{thm}

\begin{proof} The proof of theorem \ref{mainthm} consists in comparing the $T^2$
coefficient of the zeta functions associated to $f$ and $g$. Due to
the invariance theorem \ref{inv}, it is sufficient to compare
these coefficients for the simpler Nash germs given by lemma
\ref{morse}, since they are Nash equivalent and therefore blow-Nash
equivalent to $f$ and $g$ respectively. Now, the result follows from proposition \ref{propmain} below.

\end{proof}

\section{Computation of some virtual Poincar\'e polynomials}\label{partSS}

Let $X_{m,M}$ be the real algebraic subset of $\mathbb R^{m+M}$ defined by
the equation
$$\sum _{i=1}^m x_i^2 -\sum_{j=1}^M y_j^2=0.$$
In this section, we compute the value of the virtual Poincar\'e
polynomials $\beta(X_{m,M})$ in terms of $m$ and $M$. This computation is based on the degeneracy of a Leray-Serre
spectral sequence associated to the projectivisation of $X_{m,M}$.

Without lost of generality, one may assume that $m \leq M$.

\begin{prop}\label{calcul} If $m \geq 1$, then $\beta(X_{m,M})=u^{m+M-1}-u^{M-1}+u^m$.
\end{prop}

\begin{rmk}\begin{flushleft}\end{flushleft}
\begin{enumerate}
\item If $m=0$, $X_{0,M}$ is
empty, and therefore $\beta (X_{0,M})=0$. 
\item Note that if $m=1$, then the computation of $\beta
  (X_{1,M})$ is easy since $X_{1,M}$ is just a cone based on a
  sphere. Thus $\beta (X_{1,M})= 1+(u-1)(1+u^{M-1})=u^M-u^{M-1}+u$.
\item In the particular case where $m=M=2$, the computation can
  be done in a simple way using the toric structure of $X_{2,2}$. Indeed $X_{2,2}$ is isomorphic to the toric variety
  given by $XY=UV$ in $\mathbb R^4$. Therefore $X_{2,2}$ is the union of the
  orbits under the torus action, that is $X_{2,2}$ is the disjoint union of
  $(\mathbb R^*)^3$, one point, and four copies of  $(\mathbb R^*)^2$
  and $(\mathbb R^*)$. Therefore, by additivity of $\beta$,
$$\beta(X_{2,2})=(u-1)^3+4(u-1)^2+4(u-1) +1=u^3+u^2-u.$$
\end{enumerate}
\end{rmk}

The proof of proposition \ref{calcul} is based on a reduction to the
projective case. As a
preliminary step, we compute the virtual
Poincar\'e polynomial of the projective subset $Z_{m,M}$ of $\mathbb
P^{m+M-1}(\mathbb R)$ defined by the same
equation as that of $X_{m,M}$.

\begin{lemma}\label{lem1} Take $M \geq 2$. Let $Z_{m,M}$ be defined by $\sum _{i=1}^m x_i^2 -\sum_{j=1}^M
  y_j^2=0$ in $\mathbb P^{m+M-1}(\mathbb R)$. Then
$$\beta(Z_{m,M})=(1+u^{M-1})(1+u+\ldots +u^{m-1}).$$
\end{lemma}

\begin{proof} To begin with, remark that $Z_{m,M}$ is nonsingular as a real
  algebraic set. Therefore the virtual Betti numbers of $Z_{m,M}$ coincide
  with its classical Betti numbers (cf. proposition \ref{beta-nash}).

In order to compute these Betti numbers, consider the projection from $Z_{m,M}$ onto $\mathbb
  P^{m-1}(\mathbb R)$ defined by 
$$[x_1:\cdots :x_m:y_1:\cdots:y_M] \mapsto [x_1:\cdots :x_m].$$
It is well-defined since $x_1,\ldots ,x_m$ can not vanish without
cancelling $y_1,\ldots ,y_M$, and moreover it defines a fibration with
fiber isomorphic to the unit sphere $S^{M-1}$ in $\mathbb R^M$. Working with coefficients in $\mathbb Z_2$,
the cohomological Leray-Serre spectral sequence associated with this
fibration converges to the cohomology with coefficients in $\mathbb
Z_2$  of $Z_{m,M}$:
$$E_2^{p,q}=H^p(\mathbb P^{m-1}(\mathbb R), H^q (S^{M-1},\mathbb Z_2))
\Rightarrow H^{p+q}(Z_{m,M},\mathbb Z_2).$$

However $H^q (S^{M-1},\mathbb Z_2)$ is zero unless $q=0$ and $q=M-1$
for which it equals $\mathbb Z_2$. Therefore the nonzero terms (that
equals $\mathbb Z_2$) of $E_2^{p,q}$, shown in the figure below, are
localized in two lines. 

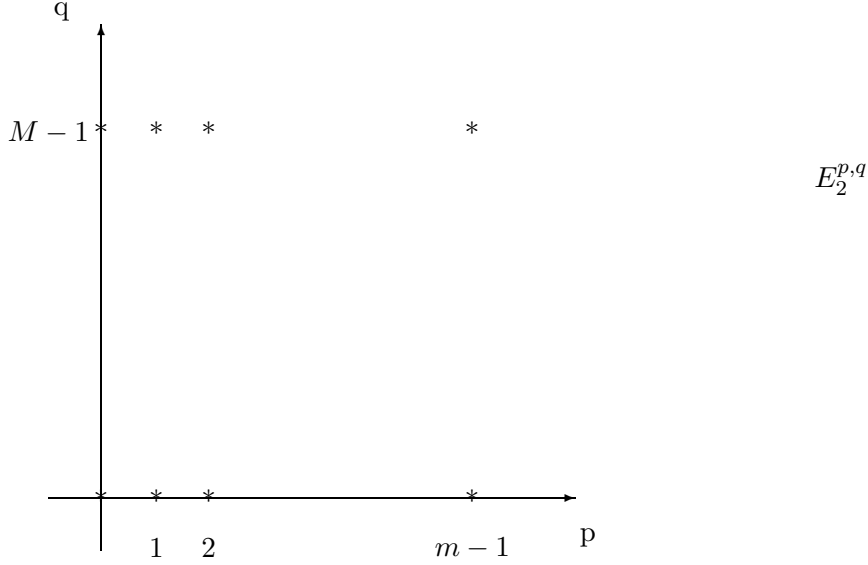
\unitlength=.7mm
\begin{picture}(150,100)(-10,10)
\put(-10,0){\vector(1,0){100}}
\put(0,-10){\vector(0,1){100}}
\put(90,-10){\makebox(5,5){p}}
\put(-10,90){\makebox(5,5){q}}
\put(0,0){\makebox(0,-1){*}}
\put(10,0){\makebox(1,-1){*}}\put(10,-10){\makebox(1,1){1}}
\put(20,0){\makebox(1,-1){*}}\put(20,-10){\makebox(1,1){2}}
\put(70,0){\makebox(1,-1){*}}
\put(0,70){\makebox(0,-1){*}}
\put(10,70){\makebox(1,-1){*}}
\put(20,70){\makebox(1,-1){*}}
\put(70,70){\makebox(1,-1){*}}
\put(-10,70){\makebox(0,-1){$M-1$}}
\put(70,-10){\makebox(1,1){$m-1$}}

\put(140,60){\makebox(1,1){$E_2^{p,q}$}}
\end{picture}

\vskip 20mm
But, as $m \leq M$ by assumption, the spectral sequence degenerates and gives the
Betti numbers of $Z_{m,M}$. More precisely, 
\begin{displaymath}
\dim H_i(Z_{m,M},\mathbb Z_2) = \left\{ \begin{array}{ll}
1 & \textrm{if~~} i\in \{0,\ldots,m-1,M-1,\ldots, m+M-2\},\\
0 & \textrm{otherwise,}
\end{array} \right.
\end{displaymath}
if
$m<M$, and in the particular case where $m=M$, then 
\begin{displaymath}
\dim H_i(Z_{m,m},\mathbb Z_2) = \left\{ \begin{array}{lll}
1 & \textrm{if~~} i\in \{0,\ldots,m-2,m,\ldots, 2m-2\},\\
2 & \textrm{if~~} i=m-1,\\
0 & \textrm{otherwise.}
\end{array} \right.
\end{displaymath}
So in general
$$\beta(Z_{m,M})=(1+u+\cdots+u^{m-1})+(u^{M-1}+\cdots+u^{m+M-2})=(1+u^{M-1})(1+u+\ldots +u^{m-1}).$$
\end{proof}

Now we explain how to compute the virtual Poincar\'e polynomial of $X_{m,M}$
in terms of that of $Z_{m,M}$.

\begin{proof}[Proof of proposition \ref{calcul}]

It suffices to notice that the projection from $X_{m,M}\setminus \{0\}$
  onto $Z_{m,M}$ is a piecewise algebraically trivial fibration with
  fiber $\mathbb R^*$. Therefore 
$$\beta(X_{m,M})=1+(u-1)\beta(Z_{m,M})$$
by additivity and multiplicativity of the virtual Poincar\'e
polynomial $\beta$.

Now, remark that $\beta(Z_{m,M})=(1+u^{M-1})\frac{u^m-1}{u-1}$, and so
$$\beta(X_{m,M})=1+(1+u^{M-1})(u^m-1)=u^{m+M-1}-u^{M-1}+u^m.$$

\end{proof}

The following corollaries, which also specify the virtual Poincar\'e
polynomial of some algebraic sets, will be usefull for computing
zeta functions with sign in section \ref{partcalculzeta}.

\begin{cor}\label{coro1} Let $X_{s,t}^1$  be the real algebraic subset of $\mathbb R^{s+t}$
  defined by the equation
$$\sum _{i=1}^s x_i^2 -\sum_{j=1}^t
y_j^2=1.$$
Assume that $s, t>0$. 
\begin{itemize}
\item If $s \leq t$, then $\beta(X_{s,t}^1)=u^{t-1}(u^s-1)$.
\item If $s > t$, then $\beta(X_{s,t}^1)=u^{t}(u^{s-1}+1)$
\end{itemize}

Moreover $\beta(X_{0,t}^1)=0$ and
$\beta(X_{s,0}^{1})=1+\cdots+u^{s-1}$ if $s\geq 1$.
\end{cor}
\begin{proof} Let us begin with the case $s \leq t$. 
Let homogenize the equation defining $X_{s,t}^{1}$. Then, we obtain a projective subset of $\mathbb
  P(\mathbb R)^{s+t}$, denoted $Z_{s,t+1}$ in lemma \ref{lem1}, whose
  affine part is isomorphic to $X_{s,t}^{1}$, and whose part at
  infinity is isomorphic to $Z_{s,t}$. Therefore
$$\beta(X_{s,t}^{1})=\beta(Z_{s,t+1})-\beta(Z_{s,t})=(1+u^t)(1+\cdots+u^{s-1})-(1+u^{t-1})(1+\cdots+u^{s-1})$$
$$=u^{t-1}(u-1)\frac{u^s-1}{u-1}=u^{t-1}(u^s-1).$$

Now, let us turn to the case $s > t$. In the same way,
 by homogenization of
  the equation defining $X_{s,t}^{1}$, we obtain a projective subset of $\mathbb
  P(\mathbb R)^{s+t}$, denoted $Z_{t+1,s}$ (and
  not $Z_{s,t+1}$ because $s \geq t+1$), whose
  affine part is isomorphic to $X_{s,t}^{1}$. Moreover, the part at
  infinity is isomorphic to $Z_{t,s}$, therefore
$$\beta(X_{s,t}^{1})=\beta(Z_{t+1,s})-\beta(Z_{t,s}),$$
and the second member can be computed thanks to lemma \ref{lem1}. More
precisely:
$$\beta(X_{s,t}^{1})=(1+u^{s-1})(1+u+\ldots +u^{t})-(1+u^{s-1})(1+u+\ldots +u^{t-1})=u^t(1+u^{s-1}).$$

Finally, remark that in the case $s=0$, then the sets considered are
either empty or isomorphic to a sphere.
\end{proof}

\begin{cor}\label{coro2} Let $X_{s,t}^{-1}$ be the real algebraic subset of $\mathbb R^{s+t}$
  defined by the equation
$$\sum _{i=1}^s x_i^2 -\sum_{j=1}^t
y_j^2=-1.$$
Assume that $s, t>0$. 
\begin{itemize}
\item If $s \geq t$, then $\beta(X_{s,t}^{-1})=u^{s-1}(u^t-1)$.
\item If $s < t$, then $\beta(X_{s,t}^{-1})=u^s(u^{t-1}+1)$.
\end{itemize}

Moreover $\beta(X_{s,0}^{-1})=0$ and
$\beta(X_{0,t}^{-1})=1+\cdots+u^{t-1}$ if $t\geq 1$.
\end{cor}

\begin{rmk} This is just a rewriting of corollary \ref{coro1}
after noticing that $X_{s,t}^{-1}=X_{t,s}^{1}$
\end{rmk}
\section{Proof of theorem \ref{mainthm}}\label{partcalculzeta}

We study the behaviour of the zeta functions of a germ of functions
$f:(\mathbb R^{d},0) \longrightarrow (\mathbb R,0)$ of
the form $$f(x_1,\ldots, x_s, y_1,\ldots,y_t,
z_1,\ldots,z_{d-s-t})=\sum_{i=1}^s x_i^2-\sum_{j=1}^t y_j^2.$$ In particular, we compute
the coefficient of $T^2$ of the naive zeta function and of the zeta
functions with sign. 

The main result, stated in proposition \ref{propmain} below, is that
the corresponding coefficients of the zeta functions with sign determine the integers
$s$ and $t$. It completes the proof of theorem \ref{mainthm}.

\begin{prop}\label{propmain} The coefficients of $T^2$ of the zeta functions with sign of a germ of analytic functions $f:(\mathbb R^{d},0) \longrightarrow (\mathbb R,0)$ of
the form $$f(x_1,\ldots, x_s, y_1,\ldots,y_t,
z_1,\ldots,z_{d-s-t})=\sum_{i=1}^s x_i^2-\sum_{j=1}^t y_j^2$$
determine $s$ and $t$.
\end{prop}

\begin{proof} The space of truncated arcs $A_2^{+1}(f)$ is isomorphic to
  $$\mathbb R^{2(d-s-t)} \times \mathbb R^{s+t}  \times
  X_{s,t}^1.$$
Actually, for an arc $(a _1 t + _1 t^2,\ldots,a_d t + b_d t^2)$ in $A_2^{+1}(f)$,
 the first term of the product corresponds to the choice of the
coefficients $a_{s+t+1}, b_{s+t+1}, \ldots, a_d,b_d$, the
second to the choice of $b_1,\ldots, b_{s+t}$ and finally $X_{s,t}^1$
to the choice of $a_1,\ldots, a_{s+t}$.

Now, putting into factor in $\beta(A_2^{+1}(f))$ the maximal power of
$u$, we remark that, because of the
possible forms of this polynomial as specified in corollary
\ref{coro1}, it remains a polynomial of the form $u^k+1$ or $u^k-1$. Now, due to corollary \ref{coro1} again, it
follows that $s=k+1$ in the former case, and $s=k$ in the latter one.

Similarly, $A_2^{-1}(f)$ is isomorphic to
  $$\mathbb R^{2(d-s-t)} \times \mathbb R^{s+t}  \times
  X_{s,t}^{-1},$$
and once more, after dividing $\beta(A_2^{-1}(f))$ by the maximal
power of $u$, we obtain a polynomial of the type $u^l+1$ or
$u^l-1$. In the former case, then $t=l+1$ whereas in the latter one $t=l$.
\end{proof}

\begin{rmk} The same method, applied to the
  coefficient of $T^2$ in the naive
  zeta function instead of the zeta functions with sign, no longer
  gives such a determination. More precisely, what we still determine is
  $m=\min \{s,t\}$ and $M=\max \{s,t\}$, unless $m=0$ or $M=m+1$ where
  we can even not specify the value of $M$ and $m,M$ respectively.
\end{rmk}

\end{document}